\RequirePackage{ifpdf}
\ifpdf % We are running pdfTeX in pdf mode
\documentclass[pdftex]{sigma}
\else
\documentclass{sigma}
\fi

\theoremstyle{definition}
\newtheorem{df}{Definition}[section]
\newtheorem{thm}[df]{Theorem}
\newtheorem{pro}[df]{Proposition}
\newtheorem{lem}[df]{Lemma}
\newtheorem{ex}[df]{Example}
\newtheorem{cor}[df]{Corollary}

\numberwithin{equation}{section}

\usepackage[all]{xy}
\input xy
\xyoption{all}

\begin{document}

\allowdisplaybreaks

\renewcommand{\thefootnote}{$\star$}

\renewcommand{\PaperNumber}{023}

\FirstPageHeading

\ShortArticleName{Epsilon Systems on Geometric Crystals of type $A_n$}

\ArticleName{Epsilon Systems on Geometric Crystals of type $\boldsymbol{A_n}$\footnote{This paper is a
contribution to the Proceedings of the Workshop ``Geometric Aspects of Discrete and Ultra-Discrete Integrable Systems'' (March 30 -- April 3, 2009, University of Glasgow, UK). The full collection is
available at
\href{http://www.emis.de/journals/SIGMA/GADUDIS2009.html}{http://www.emis.de/journals/SIGMA/GADUDIS2009.html}}}

\Author{Toshiki NAKASHIMA}

\AuthorNameForHeading{T. Nakashima}

\Address{Department of Mathematics, Sophia University, 102-8554, Chiyoda-ku, Tokyo, Japan}

\Email{\href{mailto:toshiki@mm.sophia.ac.jp}{toshiki@mm.sophia.ac.jp}}

\ArticleDates{Received September 14, 2009, in f\/inal form January 28, 2010;  Published online March 19, 2010}

\Abstract{We introduce an epsilon system on a geometric crystal of type $A_n$,
which is a~certain set of rational functions with some nice properties.
We shall show that it is equipped with
a product structure and that it is invariant
under the action of tropical R maps.}


\Classification{17B37; 17B67; 22E65; 14M15}

\renewcommand{\labelenumi}{$($\roman{enumi}$)$}
\renewcommand{\labelenumii}{$(${\rm \alph{enumii}}$)$}
\font\germ=eufm10
\newcommand{\cI}{{\mathcal I}}
\newcommand{\cA}{{\mathcal A}}
\newcommand{\cB}{{\mathcal B}}
\newcommand{\cC}{{\mathcal C}}
\newcommand{\cD}{{\mathcal D}}
\newcommand{\cE}{{\mathcal E}}
\newcommand{\cF}{{\mathcal F}}
\newcommand{\cG}{{\mathcal G}}
\newcommand{\cH}{{\mathcal H}}
\newcommand{\cJ}{{\mathcal J}}
\newcommand{\cK}{{\mathcal K}}
\newcommand{\cL}{{\mathcal L}}
\newcommand{\cM}{{\mathcal M}}
\newcommand{\cN}{{\mathcal N}}
\newcommand{\cO}{{\mathcal O}}
\newcommand{\cR}{{\mathcal R}}
\newcommand{\cP}{{\mathcal P}}
\newcommand{\cS}{{\mathcal S}}
\newcommand{\cV}{{\mathcal V}}
\newcommand{\fra}{\mathfrak a}
\newcommand{\frb}{\mathfrak b}
\newcommand{\frc}{\mathfrak c}
\newcommand{\frd}{\mathfrak d}
\newcommand{\fre}{\mathfrak e}
\newcommand{\frf}{\mathfrak f}
\newcommand{\frg}{\mathfrak g}
\newcommand{\frh}{\mathfrak h}
\newcommand{\fri}{\mathfrak i}
\newcommand{\frj}{\mathfrak j}
\newcommand{\frk}{\mathfrak k}
\newcommand{\frI}{\mathfrak I}
\newcommand{\fm}{\mathfrak m}
\newcommand{\frn}{\mathfrak n}
\newcommand{\frp}{\mathfrak p}
\newcommand{\fq}{\mathfrak q}
\newcommand{\frr}{\mathfrak r}
\newcommand{\frs}{\mathfrak s}
\newcommand{\frt}{\mathfrak t}
\newcommand{\fru}{\mathfrak u}
\newcommand{\frA}{\mathfrak A}
\newcommand{\frB}{\mathfrak B}
\newcommand{\frF}{\mathfrak F}
\newcommand{\frG}{\mathfrak G}
\newcommand{\frH}{\mathfrak H}
\newcommand{\frJ}{\mathfrak J}
\newcommand{\frN}{\mathfrak N}
\newcommand{\frP}{\mathfrak P}
\newcommand{\frT}{\mathfrak T}
\newcommand{\frU}{\mathfrak U}
\newcommand{\frV}{\mathfrak V}
\newcommand{\frX}{\mathfrak X}
\newcommand{\frY}{\mathfrak Y}
\newcommand{\frZ}{\mathfrak Z}
\newcommand{\rA}{\mathrm{A}}
\newcommand{\rC}{\mathrm{C}}
\newcommand{\rd}{\mathrm{d}}
\newcommand{\rB}{\mathrm{B}}
\newcommand{\rD}{\mathrm{D}}
\newcommand{\rE}{\mathrm{E}}
\newcommand{\rH}{\mathrm{H}}
\newcommand{\rK}{\mathrm{K}}
\newcommand{\rL}{\mathrm{L}}
\newcommand{\rM}{\mathrm{M}}
\newcommand{\rN}{\mathrm{N}}
\newcommand{\rR}{\mathrm{R}}
\newcommand{\rT}{\mathrm{T}}
\newcommand{\rZ}{\mathrm{Z}}
\newcommand{\bbA}{\mathbb A}
\newcommand{\bbB}{\mathbb B}
\newcommand{\bbC}{\mathbb C}
\newcommand{\bbG}{\mathbb G}
\newcommand{\bbF}{\mathbb F}
\newcommand{\bbH}{\mathbb H}
\newcommand{\bbP}{\mathbb P}
\newcommand{\bbN}{\mathbb N}
\newcommand{\bbQ}{\mathbb Q}
\newcommand{\bbR}{\mathbb R}
\newcommand{\bbX}{\mathbb X}
\newcommand{\bbY}{\mathbb Y}
\newcommand{\bbV}{\mathbb V}
\newcommand{\bbZ}{\mathbb Z}
\newcommand{\adj}{\operatorname{adj}}
\newcommand{\Ad}{\mathrm{Ad}}
\newcommand{\Ann}{\mathrm{Ann}}
\newcommand{\rcris}{\mathrm{cris}}
\newcommand{\ch}{\mathrm{ch}}
\newcommand{\coker}{\mathrm{coker}}
\newcommand{\diag}{\mathrm{diag}}
\newcommand{\Diff}{\mathrm{Diff}}
\newcommand{\Dist}{\mathrm{Dist}}
\newcommand{\rDR}{\mathrm{DR}}
\newcommand{\ev}{\mathrm{ev}}
\newcommand{\Ext}{\mathrm{Ext}}
\newcommand{\cExt}{\mathcal{E}xt}
\newcommand{\fin}{\mathrm{fin}}
\newcommand{\Frac}{\mathrm{Frac}}
\newcommand{\GL}{\mathrm{GL}}
\newcommand{\Hom}{\mathrm{Hom}}
\newcommand{\hd}{\mathrm{hd}}
\newcommand{\rht}{\mathrm{ht}}
\newcommand{\id}{\mathrm{id}}
\newcommand{\im}{\mathrm{im}}
\newcommand{\inc}{\mathrm{inc}}
\newcommand{\ind}{\mathrm{ind}}
\newcommand{\coind}{\mathrm{coind}}
\newcommand{\Lie}{\mathrm{Lie}}
\newcommand{\Max}{\mathrm{Max}}
\newcommand{\mult}{\mathrm{mult}}
\newcommand{\op}{\mathrm{op}}
\newcommand{\ord}{\mathrm{ord}}
\newcommand{\pt}{\mathrm{pt}}
\newcommand{\qt}{\mathrm{qt}}
\newcommand{\rad}{\mathrm{rad}}
\newcommand{\res}{\mathrm{res}}
\newcommand{\rgt}{\mathrm{rgt}}
\newcommand{\rk}{\mathrm{rk}}
\newcommand{\SL}{\mathrm{SL}}
\newcommand{\soc}{\mathrm{soc}}
\newcommand{\Spec}{\mathrm{Spec}}
\newcommand{\St}{\mathrm{St}}
\newcommand{\supp}{\mathrm{supp}}
\newcommand{\Tor}{\mathrm{Tor}}
\newcommand{\Tr}{\mathrm{Tr}}
\newcommand{\wt}{\mathrm{wt}}
\newcommand{\Ab}{\mathbf{Ab}}
\newcommand{\Alg}{\mathbf{Alg}}
\newcommand{\Grp}{\mathbf{Grp}}
\newcommand{\Mod}{\mathbf{Mod}}
\newcommand{\Sch}{\mathbf{Sch}}\newcommand{\bfmod}{{\bf mod}}
\newcommand{\Qc}{\mathbf{Qc}}
\newcommand{\Rng}{\mathbf{Rng}}
\newcommand{\Top}{\mathbf{Top}}
\newcommand{\Var}{\mathbf{Var}}
\newcommand{\gromega}{\langle\omega\rangle}
\newcommand{\lbr}{\begin{bmatrix}}
\newcommand{\rbr}{\end{bmatrix}}
\newcommand{\forb}{\bigcirc\kern-2.8ex \because}
\newcommand{\forbb}{\bigcirc\kern-3.0ex \because}
\newcommand{\forbbb}{\bigcirc\kern-3.1ex \because}
\newcommand{\cd}{commutative diagram }
\newcommand{\SpS}{spectral sequence}
\newcommand\C{\mathbb C}
\newcommand\hh{{\hat{H}}}
\newcommand\eh{{\hat{E}}}
\newcommand\F{\mathbb F}
\newcommand\fh{{\hat{F}}}
\def\ge{\frg}
\def\AA{{\mathcal A}}
\def\al{\alpha}
\def\bq{B_q(\ge)}
\def\bl{\bullet}
\def\bqm{B_q^-(\ge)}
\def\bqz{B_q^0(\ge)}
\def\bqp{B_q^+(\ge)}
\def\beneme{\begin{enumerate}}
\def\beq{\begin{equation}}
\def\beqn{\begin{eqnarray}}
\def\beqnn{\begin{eqnarray*}}
\def\bigsl{{\hbox{\fontD \char'54}}}
\def\bbra#1,#2,#3{\left\{\begin{array}{c}\hspace{-5pt}
#1;#2\\ \hspace{-5pt}#3\end{array}\hspace{-5pt}\right\}}
\def\cd{\cdots}
\def\CC{\mathbb{C}}
\def\CBL{\cB_L(\TY(B,1,n+1))}
\def\CBM{\cB_M(\TY(B,1,n+1))}
\def\CVL{\cV_L(\TY(D,1,n+1))}
\def\CVM{\cV_M(\TY(D,1,n+1))}
\def\ddd{\hbox{\germ D}}
\def\del{\delta}
\def\Del{\Delta}
\def\Delr{\Delta^{(r)}}
\def\Dell{\Delta^{(l)}}
\def\Delb{\Delta^{(b)}}
\def\Deli{\Delta^{(i)}}
\def\Delre{\Delta^{\rm re}}
\def\ei{e_i}
\def\eit{\tilde{e}_i}
\def\eneme{\end{enumerate}}
\def\ep{\epsilon}
\def\eeq{\end{equation}}
\def\eeqn{\end{eqnarray}}
\def\eeqnn{\end{eqnarray*}}
\def\fit{\tilde{f}_i}
\def\FF{{\rm F}}
\def\ft{\tilde{f}}
\def\gau#1,#2{\left[\begin{array}{c}\hspace{-5pt}#1\\
\hspace{-5pt}#2\end{array}\hspace{-5pt}\right]}
\def\gl{\hbox{\germ gl}}
\def\hom{{\hbox{Hom}}}
\def\ify{\infty}
\def\io{\iota}
\def\kp{k^{(+)}}
\def\km{k^{(-)}}
\def\llra{\relbar\joinrel\relbar\joinrel\relbar\joinrel\rightarrow}
\def\lan{\langle}
\def\lar{\longrightarrow}
\def\max{{\rm max}}
\def\lm{\lambda}
\def\Lm{\Lambda}
\def\mapright#1{\smash{\mathop{\longrightarrow}\limits^{#1}}}
\def\Mapright#1{\smash{\mathop{\Longrightarrow}\limits^{#1}}}
\def\mm{{\bf{\rm m}}}
\def\nd{\noindent}
\def\nn{\nonumber}
\def\nnn{\hbox{\germ n}}
\def\catob{{\mathcal O}(B)}
\def\oint{{\mathcal O}_{\rm int}(\ge)}
\def\ot{\otimes}
\def\op{\oplus}
\def\opi{\ovl\pi_{\lm}}
\def\osigma{\ovl\sigma}
\def\ovl{\overline}
\def\plm{\Psi^{(\lm)}_{\io}}
\def\qq{\qquad}
\def\q{\quad}
\def\QQ{\mathbb Q}
\def\qi{q_i}
\def\qii{q_i^{-1}}
\def\ra{\rightarrow}
\def\ran{\rangle}
\def\rlm{r_{\lm}}
\def\ssl{\hbox{\germ sl}}
\def\slh{\widehat{\ssl_2}}
\def\ti{t_i}
\def\tii{t_i^{-1}}
\def\til{\tilde}
\def\tm{\times}
\def\ttt{\frt}
\def\TY(#1,#2,#3){#1^{(#2)}_{#3}}
\def\ua{U_{\AA}}
\def\ue{U_{\vep}}
\def\uq{U_q(\ge)}
\def\uqp{U'_q(\ge)}
\def\ufin{U^{\rm fin}_{\vep}}
\def\ufinp{(U^{\rm fin}_{\vep})^+}
\def\ufinm{(U^{\rm fin}_{\vep})^-}
\def\ufinz{(U^{\rm fin}_{\vep})^0}
\def\uqm{U^-_q(\ge)}
\def\uqmq{{U^-_q(\ge)}_{\bf Q}}
\def\uqpm{U^{\pm}_q(\ge)}
\def\uqq{U_{\bf Q}^-(\ge)}
\def\uqz{U^-_{\bf Z}(\ge)}
\def\ures{U^{\rm res}_{\AA}}
\def\urese{U^{\rm res}_{\vep}}
\def\uresez{U^{\rm res}_{\vep,\ZZ}}
\def\util{\widetilde\uq}
\def\uup{U^{\geq}}
\def\ulow{U^{\leq}}
\def\bup{B^{\geq}}
\def\blow{\ovl B^{\leq}}
\def\vep{\varepsilon}
\def\vp{\varphi}
\def\vpi{\varphi^{-1}}
\def\VV{{\mathcal V}}
\def\xii{\xi^{(i)}}
\def\Xiioi{\Xi_{\io}^{(i)}}
\def\W1{W(\varpi_1)}
\def\WW{{\mathcal W}}
\def\wt{{\rm wt}}
\def\wtil{\widetilde}
\def\what{\widehat}
\def\wpi{\widehat\pi_{\lm}}
\def\ZZ{\mathbb Z}

\def\m@th{\mathsurround=0pt}
\def\fsquare(#1,#2){
\hbox{\vrule$\hskip-0.4pt\vcenter to #1{\normalbaselines\m@th
\hrule\vfil\hbox to #1{\hfill$\scriptstyle #2$\hfill}\vfil\hrule}$\hskip-0.4pt
\vrule}}

\newcommand{\cmt}{\marginpar}
\newcommand{\seteq}{\mathbin{:=}}
\newcommand{\cl}{\colon}
\newcommand{\be}{\begin{enumerate}}
\newcommand{\ee}{\end{enumerate}}
\newcommand{\bnum}{\be[{\rm (i)}]}
\newcommand{\enum}{\ee}
\newcommand{\ro}{{\rm(}}
\newcommand{\rf}{{\rm)}}
\newcommand{\set}[2]{\left\{#1\,\vert\,#2\right\}}
\newcommand{\sbigoplus}{{\mbox{\small{$\bigoplus$}}}}
\newcommand{\ba}{\begin{array}}
\newcommand{\ea}{\end{array}}
\newcommand{\on}{\operatorname}
\newcommand{\eq}{\begin{eqnarray}}
\newcommand{\eneq}{\end{eqnarray}}
\newcommand{\hs}{\hspace*}

\section{Introduction}

In the theory of crystal bases, the piecewise-linear functions
$\vep_i$ and $\vp_i$ play many crucial roles, e.g.,
description of highest weight vectors, tensor product of crystals,
extremal vectors, etc.
There exist counterparts for geometric crystals,
denoted also $\{\vep_i\}$, which are rational functions with
several nice properties, indeed,
they are needed to describe the product structure of geo\-met\-ric crystals
(see Section~\ref{section2}) and in $\ssl_2$-case, the universal tropical R map
is presented by using them~\cite{N5}.

In \cite{BK}, higher objects $\vep_{i,j}$ and $\vep_{j,i}$ are
introduced in order to prove the existence of product structure
of geometric crystals induced from unipotent crystals, which satisfy
the relation
$\vep_i\vep_j=\vep_{i,j}+\vep_{j,i}$ if the vertices $i$ and $j$ are
simply laced.
It motivates us to def\/ine further higher objects,
``epsilon system''.

The aim of the article is to def\/ine an ``epsilon system'' for type $A_n$
and reveal its basic properties, e.g., product structures and
invariance under the action of tropical R maps.
An epsilon system is a certain set of rational functions on a
geometric crystal, which satisfy some relations with each other
and have simple forms of the action by $e_i^c$'s.

We found its prototype on the geometric crystal of
the opposite Borel subgroup $B^-\subset SL_{n+1}(\bbC)$.
In that case, indeed, the epsilon system is realized as a set of
matrix elements and minor determinants of unipotent part of
a group element in $B^-$ (see Section~\ref{borel}).
Therefore, we know that geometric crystals induced from unipotent crystals
are equipped with an epsilon system naturally.

We shall introduce two remarkable properties of epsilon system:
One is a product structure of epsilon systems.
That is, for two geometric crystals with epsilon systems,
say $\bbX$ and $\bbY$, there exists canonically
an epsilon system on the product of geometric crystals
$\bbX\times\bbY$ (see Section~\ref{subsec-prod}).
The other is an invariance by tropical R maps:
Let $\cR:\bbX\times\bbY\to\bbY\times\bbX$ be a tropical R map
(see Section~\ref{trop-sec}) and $\vep^{\bbX\times\bbY}_J$
(resp.\ $\vep^{\bbY\times\bbX}_J$)
an arbitrary element in the epsilon system on $\bbX\times\bbY$
(resp. $\bbY\times\bbX$) obtained from the ones on $\bbX$ and $\bbY$.
Then, we have the invariant property:
\[
 \vep^{\bbY\times\bbX}_J(\cR(x,y))=\vep^{\bbX\times\bbY}_J(x,y).
\]
In the last section, we shall give an application of these results,
which shows the uniqueness of tropical R map on some geometric crystals
of type $A^{(1)}_n$.

Since the $\ssl_2$-universal tropical R map is presented by using
the rational functions $\{\vep_i\}$ as mentioned above, we expect
that epsilon systems would be a key to f\/ind universal tropical R maps
of higher ranks.
For further aim, we would like to extend this notion to other
simple Lie algebras, e.g., $B_n$, $C_n,$ and $D_n$.
These problems will be discussed elsewhere.

\section{Geometric crystals and unipotent crystals}\label{section2}

The notations and def\/initions here follow
\cite{KNO, Kac,KP,PK, Ku, N, N2, N3}.

\subsection{Geometric crystals}
\label{KM}
Fix a symmetrizable generalized Cartan matrix
 $A=(a_{ij})_{i,j\in I}$ with a f\/inite index set $I$.
Let $(\ttt,\{\al_i\}_{i\in I},\{h_i\}_{i\in I})$
be the associated
root data
satisfying $\al_j(h_i)=a_{ij}$.
Let $\ge=\ge(A)=\lan \ttt,e_i,f_i(i\in I)\ran$
be the Kac--Moody Lie algebra associated with $A$
\cite{Kac}.
Let $P\subset \ttt^*$ (resp.\ $Q:=\oplus_i\bbZ\al_i$,
$Q^\vee:=\oplus_i\bbZ h_i$) be a weight (resp.\ root, coroot)
lattice such that
$\bbC\ot P=\ttt^*$ and
$P\subset
\{\lm\,|\, \lm(Q^\vee)\subset \bbZ\}$,
whose element is called a weight.

Def\/ine the simple ref\/lections $s_i\in{\rm Aut}(\ttt)$ $(i\in I)$ by
$s_i(h)\seteq h-\al_i(h)h_i$, which generate the Weyl group $W$.
Let $G$ be the Kac--Moody group associated
with $(\ge,P)$ \cite{KP,PK}.
Let $U_{\al}\seteq\exp\ge_{\al}$ $(\al\in \Delre)$
be the one-parameter subgroup of $G$.
The group $G$ (resp.\ $U^\pm$) is generated by
$\{U_{\al}|\al\in \Delre\}$
(resp.\ $\{U_\al|\al\in\Delre\cap(\oplus_i\pm\bbZ\al_i)$).
Here $U^\pm$ is a unipotent subgroup of $G$.
For any $i\in I$, there exists
a unique group homomorphism
$\phi_i\cl SL_2(\bbC)\rightarrow G$ such that
\[
\phi_i\left(
\left(
\begin{array}{cc}
1&t\\
0&1
\end{array}
\right)\right)=\exp(t e_i),\qquad
 \phi_i\left(
\left(
\begin{array}{cc}
1&0\\
t&1
\end{array}
\right)\right)=\exp(t f_i), \qquad t\in\bbC.
\]
Set $\al^\vee_i(c)\seteq
\phi_i\left(\left(
\begin{smallmatrix}
c&0\\
0&c^{-1}\end{smallmatrix}\right)\right)$,
$x_i(t)\seteq\exp{(t e_i)}$, $y_i(t)\seteq\exp{(t f_i)}$,
$G_i\seteq\phi_i(SL_2(\bbC))$,
$T_i\seteq \alpha_i^\vee(\bbC^\times)$
and
$N_i\seteq N_{G_i}(T_i)$. Let
$T$ be the subgroup of $G$
with $P$ as its weight lattice
which is called a~{\it maximal torus} in $G$, and let
$B^{\pm}(\supset T)$ be the Borel subgroup of $G$.
We have the isomorphism
$\phi:W\mapright{\sim}N/T$ def\/ined by $\phi(s_i)=N_iT/T$.
An element $\ovl s_i:=x_i(-1)y_i(1)x_i(-1)$ is in
$N_G(T)$, which is a representative of
$s_i\in W=N_G(T)/T$.

\begin{df}
\label{def-gc}
Let $X$ be an ind-variety over $\bbC$, $\gamma_i$ and $\vep_i$
$(i\in I)$ rational functions on $X$, and
$e_i:\bbC^\times\times X\to X$ a rational $\bbC^\times$-action.
A quadruple $(X,\{e_i\}_{i\in I},\{\gamma_i,\}_{i\in I},
\{\vep_i\}_{i\in I})$ is a
$G$ (or~$\ge$)-{\it geometric crystal}
if
\begin{enumerate}\itemsep=0pt
%\bnum
\item
$(\{1\}\times X)\cap {\rm dom}(e_i)$
is open dense in $\{1\}\times X$ for any $i\in I$, where
${\rm dom}(e_i)$ is the domain of def\/inition of
$e_i\cl\C^\times\times X\to X$.
\cmt{changed}
\item
The rational functions  $\{\gamma_i\}_{i\in I}$ satisfy
$\gamma_j(e^c_i(x))=c^{a_{ij}}\gamma_j(x)$ for any $i,j\in I$.
\item
$e_i$ and $e_j$ satisfy the following relations:
 \begin{alignat*}{3}
& e^{c_1}_{i}e^{c_2}_{j}
=e^{c_2}_{j}e^{c_1}_{i} \qquad &&
{\rm if }\ a_{ij}=a_{ji}=0,& \\
&  e^{c_1}_{i}e^{c_1c_2}_{j}e^{c_2}_{i}
=e^{c_2}_{j}e^{c_1c_2}_{i}e^{c_1}_{j} \qquad &&
{\rm if }\ a_{ij}=a_{ji}=-1, & \\
& e^{c_1}_{i}e^{c^2_1c_2}_{j}e^{c_1c_2}_{i}e^{c_2}_{j}
=e^{c_2}_{j}e^{c_1c_2}_{i}e^{c^2_1c_2}_{j}e^{c_1}_{i}\qquad &&
{\rm if }\ a_{ij}=-2,\ a_{ji}=-1, & \\
&e^{c_1}_{i}e^{c^3_1c_2}_{j}e^{c^2_1c_2}_{i}
e^{c^3_1c^2_2}_{j}e^{c_1c_2}_{i}e^{c_2}_{j}
=e^{c_2}_{j}e^{c_1c_2}_{i}e^{c^3_1c^2_2}_{j}e^{c^2_1c_2}_{i}
e^{c^3_1c_2}_je^{c_1}_i\qquad &&
{\rm if }\ a_{ij}=-3,\ a_{ji}=-1. &
\end{alignat*}
\item
The rational functions $\{\vep_i\}_{i\in I}$ satisfy
$\vep_i(e_i^c(x))=c^{-1}\vep_i(x)$ and
$\vep_i(e_j^c(x))=\vep_i(x)$ if $a_{i,j}=a_{j,i}=0$.
\end{enumerate}
\end{df}

The relations in (iii) is called
{\it Verma relations}.
If $\chi=(X,\{e_i\},\{\gamma_i\},\{\vep_i\})$
satisf\/ies the conditions (i), (ii) and (iv),
we call $\chi$ a {\it pre-geometric crystal}.

\medskip

\noindent
{\bf Remark.}
The last condition (iv) is slightly modif\/ied from
\cite{KNO, N, N2, N3, N5} since all $\vep_i$ appearing in these references
satisfy the new condition `$\vep_i(e_j^c(x))=\vep_i(x)$ if
$a_{i,j}=a_{j,i}=0$'
and we need this condition
to def\/ine ``epsilon systems'' later.

\subsection{Unipotent crystals}

In the sequel, we denote the unipotent subgroup
$U^+$ by $U$.
We def\/ine unipotent crystals (see \mbox{\cite{BK,N}})
associated to Kac--Moody groups.
\begin{df}
Let $X$ be an ind-variety over $\bbC$ and
$\al:U\times X\rightarrow X$ be a rational $U$-action
such that $\al$ is def\/ined on $\{e\}\times X$. Then,
the pair ${\bf X}=(X,\al)$ is called a $U$-{\it variety}.
For $U$-varieties ${\bf X}=(X,\al_X)$
and ${\bf Y}=(Y,\al_Y)$,
a rational map
$f:X\rightarrow Y$ is called a
$U$-{\it morphism} if it commutes with
the action of $U$.
\end{df}
Now, we def\/ine a $U$-variety structure on $B^-=U^-T$.
As in \cite{Ku},  the Borel subgroup $B^-$ is an ind-subgroup of $G$ and hence
an ind-variety over $\bbC$.
The multiplication map in $G$ induces the open embedding;
$ B^-\times U\hookrightarrow G,$
which is a birational map.
Let us denote the inverse birational map by
$ g:G\longrightarrow B^-\times U$
and let rational maps
$\pi^-:G\rightarrow B^-$ and
$\pi:G\rightarrow U$ be
$\pi^-:={\rm proj}_{B^-}\circ g$
and $\pi:={\rm proj}_U\circ g$.
Now we def\/ine the rational $U$-action $\al_{B^-}$ on $B^-$ by
\[
 \al_{B^-}:=\pi^-\circ m: \ U\times B^-\longrightarrow B^-,
\]
where $m$ is the multiplication map in $G$.
Then we get $U$-variety ${\bf B}^-=(B^-,\al_{B^-})$.
\begin{df}\label{uni-def}\qquad
\begin{enumerate}\itemsep=0pt
\item
Let ${\bf X}=(X,\al)$
be a $U$-variety and $f:X \rightarrow B^-$
a $U$-morphism.
The pair $({\bf X}, f)$ is called
a {\it unipotent $G$-crystal}
or, for short, {\it unipotent crystal}.
\item
Let $({\bf X},f_X)$ and $({\bf Y},f_Y)$
be unipotent crystals.
A $U$-morphism $g:\bf X\ra \bf Y$ is called a~{\it morphism of
unipotent crystals} if $f_X=f_Y\circ g$.
In particular, if $g$ is a birational map
of ind-varieties, it is called an {\it isomorphism of
unipotent crystals}.
\end{enumerate}
\end{df}

We def\/ine a product of
unipotent crystals following \cite{BK}.
For unipotent crystals $({\bf X},f_X)$, $({\bf Y},f_Y)$,
def\/ine a rational map
$\al_{X\times Y}:U\tm X\tm Y\rightarrow X\tm Y$ by
\begin{equation*}
\al_{X\tm Y}(u,x,y):=(\al_X(u,x),\al_Y(\pi(u\cdot f_X(x)),y)).
%\label{XY}
\end{equation*}

\begin{thm}[\cite{BK}] \qquad \it
\label{prod}
\begin{enumerate}\itemsep=0pt
\item
The rational map $\al_{X\tm Y}$ defined above
is a rational $U$-action
on $X\tm Y$.
\item
Let ${\bf m}:B^-\tm B^-\rightarrow B^-$
be a multiplication map
and $f=f_{X\tm Y}:X\tm Y\rightarrow B^-$ be the
rational map defined by
\[
f_{X\tm Y}:={\bf m}\circ( f_X\tm f_Y).
\]
Then $ f_{X\tm Y}$ is a $U$-morphism and
$({\bf X\tm Y}, f_{X\tm Y})$ is a~unipotent crystal,
which we call a~product of unipotent crystals
$({\bf X},f_X)$ and $({\bf Y},f_Y)$.
\item
Product of unipotent crystals is associative.
\end{enumerate}
\end{thm}

\subsection{From unipotent crystals to geometric crystals}

For $i\in I$,
set $U^\pm_i:=U^\pm\cap \bar s_i U^\mp\bar s_i^{-1}$ and
$U_\pm^i:=U^\pm\cap \bar s_i U^\pm\bar s_i^{-1}$.
Indeed, $U^\pm_i=U_{\pm \al_i}$.
Set
\[
 Y_{\pm\al_i}:=
\lan x_{\pm\al_i}(t)U_{\al}x_{\pm\al_i}(-t)
|t\in\bbC,\,\,\al\in \Delta^{\rm re}_{\pm}\setminus
\{\pm\al_i\}\ran,
\]
where $x_{\al_i}(t):=x_i(t)$ and $x_{-\al_i}(t):=y_i(t)$.
We have the unique
decomposition;
\[
 U^-=U_i^-\cdot Y_{\pm\al_i}=U_{-\al_i}\cdot U^i_-.
\]
By using this decomposition, we get
the canonical projection
$\xi_i:U^-\rightarrow U_{-\al_i}$.
Now, we def\/ine the function on $U^-$ by
\[
\chi_i:=y_i^{-1}\circ\xi_i: \
U^-\longrightarrow U_{-\al_i}\mapright{\sim}
\bbC,
\]
and extend this to the function on
$B^-$ by $\chi_i(u\cdot t):=\chi_i(u)$ for
$u\in U^-$ and $t\in T$.
For a~unipotent $G$-crystal $({\bf X},f_X)$, we def\/ine a function
$\vep_i:=\vep_i^X:X\rightarrow \bbC$ by
\[
\vep_i:=\chi_i\circ{f_X},
%\label{chii}
\]
and a rational function
$\gamma_i:X\rightarrow \bbC$ by
\begin{equation*}
\gamma_i:=
\al_i\circ{\rm proj}_T\circ{ f_X}: \ X\rightarrow B^-\rightarrow T\to\bbC,
%\label{gammax}
\end{equation*}
where ${\rm proj}_T$ is the canonical projection.

\medskip

\noindent
{\bf Remark.}
Note that the function $\vep_i$ is denoted by $\vp_i$ in
\cite{BK, N}.

\medskip

Suppose that the function $\vep_i$ is not identically zero on $X$.
We def\/ine a morphism $e_i$: \mbox{$\bbC^\tm  \tm  X \rightarrow X$} by
\begin{equation*}
e^c_i(x):=x_i
\left({\frac{c-1}{\vep_i(x)}}\right)(x).
%\label{ei}
\end{equation*}

\begin{thm}[\cite{BK}]\label{U-G}\it
For a unipotent $G$-crystal $({\bf X},f_X)$,
suppose that
the function $\vep_i$ is not identically zero
for any $i\in I$.
Then the rational functions $\gamma_i,\vep_i:X\rightarrow \bbC$
and
$e_i:\bbC^\tm\tm  X\rightarrow X$ as above
define a geometric
$G$-crystal $({\bf X},\{e_i\}_{i\in I},\{\gamma_i\}_{i\in I},
\{\vep_i\}_{i\in I})$,
which is called the induced geometric $G$-crystals by
unipotent $G$-crystal $({\bf X},f_X)$.
\end{thm}

\begin{pro}\it
For unipotent $G$-crystals $({\bf X},f_X)$
and $({\bf Y},f_Y)$, set
the product $({\bf Z},f_Z):=({\bf X},f_X)
\tm({\bf Y},f_Y)$, where
$Z=X\tm Y$. Let $(Z,\{e^Z_i\}_{i\in I},\{\gamma^Z_i\}_{i\in I},
\{\vep^Z_i\}_{i\in I})$
be the induced geometric
$G$-crystal from $({\bf Z},f_Z)$.
Then we obtain:
\begin{enumerate}\itemsep=0pt
\item For each $i\in I$, $(x,y)\in Z$,
\begin{equation*}
\gamma^Z_i(x,y)=\gamma^X_i(x)\gamma^Y_i(y),\qquad
\vep^Z_i(x,y)=
\vep^X_i(x)+\frac{\vep^Y_i(y)}{\gamma^X_i(x)}.
%\label{gamma-zxy}
\end{equation*}
\item
For any $i\in I$, the action
$e^Z_i:\bbC^\tm\tm Z\rightarrow Z$ is given
by: $(e^Z_i)^c(x,y)=((e^X_i)^{c_1}(x),(e^Y_i)^{c_2}(y))$, where
\begin{equation*}
c_1=
\frac{c\gamma^X_i(x)\vep^X_i(x)+\vep^Y_i(y)}
{\gamma^X_i(x)\vep^X_i(x)+\vep^Y_i(y)},\qquad
c_2=
\frac{c(\gamma^X_i(x)\vep^X_i(x)+\vep^Y_i(y))}
{c\gamma^X_i(x)\vep^X_i(x)+\vep^Y_i(y)}.
%\label{c1c2}
\end{equation*}
\end{enumerate}
\end{pro}
Here note that $c_1c_2=c$.
The formula $c_1$ and $c_2$ in \cite{BK}
seem to be dif\/ferent from ours.

\section{Prehomogeneous geometric crystal}

\begin{df} Let $\chi=(X,\{e^c_i\},
\{\gamma_i\},\{\vep_i\})$ be
a geometric crystal. We say that
$\chi$ is {\it prehomogeneous}
if there exists
a Zariski open dense subset $\Omega\subset X$
which is an orbit by the actions of the
$e_i^c$'s.
\end{df}

\begin{lem}[\cite{KNO2}]\label{uniq}\it
Let $\chi_j=(X_j,\{e^c_i\},
\{\gamma_i\},\{\vep_i\})$ $(j=1,2)$
be  prehomogeneous geometric crystals.
Let
$\Omega_1\subset X_1$ be an open dense
orbit in $X_1$.
For isomorphisms of geometric
crystals $\phi, \phi'\cl\chi_1\to\chi_2$,
suppose that there exists $p_1\in \Omega_1$
such that
$\phi(p_1)=\phi'(p_1)\in X_2$.
Then, we have $\phi=\phi'$ as rational maps.
\end{lem}

\begin{thm}[\cite{KNO2}]\label{conn-preh}\it
Let $\chi=(X,\{e^c_i\},
\{\gamma_i\},\{\vep_i\})$ be
a finite-dimensional positive
geometric crystal with
the positive structure $\theta\cl (\bbC^\tm)^{\dim(X)}\to X$ and
$B\seteq
 {\mathcal UD}_\theta(\chi)$  the
crystal obtained as the ultra-discretization
of $\chi$.
If $B$ is a connected crystal,
then $\chi$ is prehomogeneous.
\end{thm}

In \cite{KNO,KNO2}, we showed that
ultra-discretization  of the af\/f\/ine geometric crystal
$\cV(\ge)_{l}$ $(l>0)$
is a limit of perfect crystal $B_\infty(\ge^L)$,
where $\ge^L$ is the Langlands dual of $\ge$. Since for any $k\in\bbZ_{>0}$
a tensor product $B_\ify(\ge^L)^{\ot k}$ is connected by the perfectness of
$B_\ify(\ge^L)$ and we have the isomorphism of crystals
\[
{\mathcal UD}(\cV(\ge)_{L_1}\times\cd
\times\cV(\ge)_{L_k})
\cong B_{\ify}\big(\ge^L\big)^{\ot k},
\qquad
L_1,\dots,L_k>0,
\]
by Theorem \ref{conn-preh} we obtain the following:

\begin{cor}[\cite{KNO2}]\label{preh}\it
$\cV(\ge)_{L_1}\times\cd
\times\cV(\ge)_{L_k}$ is prehomogeneous.
\end{cor}

%%%%%%%%%%%%%%%%%%%%%%%%%%%%%%%%%%%%%%%%%%%%%%%%%%%%%%%%
\renewcommand{\thesection}{\arabic{section}}
\section{Tropical R maps}\label{trop-sec}
\setcounter{equation}{0}
\renewcommand{\theequation}{\thesection.\arabic{equation}}

\begin{df}
\label{def-trop-r}
Let $\{X_{\lm}\}_{\lm\in\Lm}$ be a family of geometric crystals
with the product structures, where $\Lm$ is an index set.
A birational map
${\cal R_{\lm\mu}}:X_\lm\times X_\mu\longrightarrow
X_\mu\times X_\lm$ is called a
{\it tropical R map}  if they satisfy:
\begin{gather}
\big(e^{X_\mu\times X_\lm}_i\big)^c\circ{\cal R}_{\lm\mu}
={\cal R}_{\lm\mu}\circ \big(e^{X_\lm\times X_\mu}_i\big)^c,
\label{r1}\\
\vep^{X_\lm\times X_\mu}_i
=\vep^{X_\mu\times X_\lm}_i\circ{\cal R}_{\lm\mu},
\label{r2}\\
\gamma_i^{X_\lm\times X_\mu}
=\gamma_i^{X_\mu\times X_\lm}
\circ{\cal R}_{\lm\mu},\label{r3}\\
{\cal R}_{\lm\mu}{\cal R}_{\mu\nu}
{\cal R}_{\lm\mu}
={\cal R}_{\mu\nu}{\cal R}_{\lm\mu}
{\cal R}_{\mu\nu}
\label{r4}
\end{gather}
for any $i\in I$ and $\lm,\mu,\nu\in\Lm$.
\end{df}

Tropical R maps for certain
af\/f\/ine geometric crystals
of type $\TY(A,1,n)$, $\TY(B,1,n)$, $\TY(D,1,n)$,
$\TY(D,2,n+1)$, $\TY(A,2,2n-1)$ and $\TY(A,2,2n)$
are described explicitly~\cite{KNO2,KOTY}.

The following is immediate from Lemma \ref{uniq} and
Corollary \ref{preh}.
\begin{thm}[\cite{KNO2}]\label{uniqueness}\it
Let $\cR, \cR':\cV_L\times\cV_M\to \cV_M\times\cV_L$
be tropical R maps.
Suppose that there exists $p\in \cV_L\times \cV_M$ such that
$\cR(p)=\cR'(p)$. Then we have
$\cR=\cR'$ as birational maps.
\end{thm}

Let us introduce an example of a tropical R map of
type $\TY(A,1,n)$.
\begin{ex}\rm
\label{ex-A}
Set $\cB_L:=\{l=(l_1,\dots,l_{n+1})\,|\,l_1l_2\cd l_{n+1}=L\}$,
which is equipped with an
$\TY(A,1,n)$-geometric crystal structure by:
\[
e_i^c(l)=(\dots,c\,l_i,c^{-1}\,l_{i+1},\dots),\qquad
\gamma_i(l)=l_i/l_{i+1},\qquad \vep_i(l)=l_{i+1},\qquad
i=0,1,\dots,n.
 \]
The tropical R map on $\{\cB_L\}_{L\in\bbR_{>0}}$
is given by {\rm \cite{KOTY}}:
\begin{gather}
 \cR: \ \cB_L\times\cB_M\to\cB_M\times\cB_L, \qquad \cR(l,m)=(l',m'),\label{tropr-A}\\
l'_i:=m_i\frac{P_i(l,m)}{P_{i-1}(l,m)},\qquad
m'_i:=l_i\frac{P_{i-1}(l,m)}{P_{i}(l,m)},\qquad
\text{where}\quad
P_i(l,m):=\sum_{k=1}^{n+1}
\prod_{j=k}^{n+1}l_{i+j}\prod_{j=1}^{k}m_{i+j}.\nonumber
\end{gather}
\end{ex}

\noindent
{\bf Remark.}
In the case $\ge=A^{(1)}_n$, we have $\cB_L\cong \cV_L$.

\section{Epsilon systems}\label{ep-sys}

\subsection{Def\/inition of epsilon systems}
\begin{df}
For an interval
$J:=\{s,s+1,\dots,t-1,t\}\subset I=\{1,2,\dots,n\}$
a set of intervals $P=\{I_1,\dots,I_k\}$ is called
a {\it partition} of $J$ if disjoint intervals
$I_1,\dots,I_k$
satisfy $I_1\sqcup\cd\sqcup I_k=J$ and
max($I_j$)$<$min($I_{j+1}$) $(j=1,\dots,k-1)$.
\end{df}

For a partition $P=\{I_1,\dots,I_k\}$ of some interval $J$,
set $l(P):=k$  and called the {\it length} of $P$.
Let $\cal J$ be the set of all intervals in $I$.
For an interval $J\in\cal J$, def\/ine
\[
\cP_J:=\{P\, |\, P\text{ is a partition of }J\}.
\]
For a partition $P=\{I_1,\dots,I_k\}$ and
symbols $\vep_{I_j}$ ($j=1,\dots,k$), def\/ine
\[
 \vep_P:=\vep_{I_1}\cdot\vep_{I_2}\cdots \vep_{I_k}.
\]

\begin{df}
\label{def-ep}
For an $A_n$-geometric crystal ${\bbX}
=(X,\{e_i^c\},\{\gamma_i\},\{\vep_i\})$, the set of
rational functions on $X$, say
$E=\{\vep_J, \vep^*_J|J\in\cal J\text{ is an interval }\}$, is
called an {\it epsilon system} of $\bbX$ if they
satisfy the following:
\begin{gather}
 \vep_J(e_i^c x)=\begin{cases}c^{-1}\vep_J(x)
&\text{if} \ i=s,\\
\vep_J(x)&\text{if} \ i\ne s-1,s,t+1,
\end{cases}\nonumber\\
\vep^*_J(e_i^c x)=\begin{cases}c^{-1}\vep^*_J(x)
&\text{if} \ i=t,\\
\vep^*_J(x)&\text{if} \ i\ne s-1,t,t+1,
\end{cases}
\label{ep-ei}\\
\vep^*_J=\sum_{P\in \cP_J}(-1)^{|J|-l(P)}\vep_P,\qquad
\text{for any interval} \ J=\{s,s+1,\dots,t\}\subset I,
\label{ep*ep}
\end{gather}
and we set $\vep_J:=\vep_i$ for $J=\{i\}$,
which is originally equipped with $\bbX$. Note that $\vep^*_i=\vep_i$.
We call a geometric crystal with an epsilon system an
$\vep$-{\it geometric crystal}.
\end{df}

The actions of $e^c_{s-1}$ and $e^c_{t+1}$ will be
described explicitly below, which are derived from (\ref{ep-ei})
and (\ref{ep*ep}).

\begin{pro}\label{well-def}\it
The above definition is well-defined, that is,
for any $J\in{\cal J}$ we have
\begin{gather}
 \vep_J(e^{c_1}_{i}e^{c_2}_{j}x)
=\vep_J(e^{c_2}_{j}e^{c_1}_{i}x),\qquad
\vep^*_J(e^{c_1}_{i}e^{c_2}_{j}x)
=\vep^*_J(e^{c_2}_{j}e^{c_1}_{i}x)\qquad
\text{if}\quad \,a_{ij}=a_{ji}=0,
\label{well-0}\\
\vep_J(e^{c_1}_{i}e^{c_1c_2}_{j}e^{c_2}_{i}x)
=\vep_J(e^{c_2}_{j}e^{c_1c_2}_{i}e^{c_1}_{j}x),\nonumber\\
\vep^*_J(e^{c_1}_{i}e^{c_1c_2}_{j}e^{c_2}_{i}x)
=\vep^*_J(e^{c_2}_{j}e^{c_1c_2}_{i}e^{c_1}_{j}x),
\qquad
\text{if}\quad a_{ij}=a_{ji}=-1.
\label{well-1}
\end{gather}
More precisely, we claim that if we calculate the both sides of
the above equations by using \eqref{ep-ei} and \eqref{ep*ep},
they coincide with each other.
\end{pro}

The proof will be given in the next subsection.
\begin{ex}\label{ex-1}\qquad \rm
\begin{enumerate}\itemsep=0pt
\item $I=\{1,2\}$: $E=\{\vep_1,\vep_2,\vep_{12},
\vep_{21}=\vep^*_{12}\}$ with the relation
$\vep_{21}=\vep_1\vep_2-\vep_{12}$.
\item
$I=\{1,2,3\}$: $E=\{\vep_1,\vep_2,\vep_3,\vep_{12},
\vep^*_{12},\vep_{23},\vep^*_{23},\vep_{123},\vep^*_{123}\}$
with the relations for rank 2 and
$\vep^*_{123}=\vep_{123}-\vep_1\vep_{23}-\vep_{12}\vep_3+
\vep_1\vep_2\vep_3$.
\end{enumerate}
\end{ex}

The following lemma will be needed in the sequel.
Set
$[l,m]:=\{l,l+1,\dots,m-1,m\}$ $(l\leq m)$ and
if $l>m$, $[l,m]:=\varnothing$. We also put $\vep_\varnothing(x)=1$.
Let $\{\vep_J\}_{J\in \cJ}$ be an epsilon system on
an $\vep$-geometric crystal $\bbX$.
\begin{lem}\label{equiv}\it
One of the following relations \eqref{non-star} and
\eqref{star-non} is equivalent to \eqref{ep*ep}:
\begin{gather}
\sum_{j=s-1}^t(-1)^{j}\vep_{[s,j]}\vep^*_{[j+1,t]}=0,
\label{non-star}\\
\sum_{j=s-1}^t(-1)^{j}\vep^*_{[s,j]}\vep_{[j+1,t]}=0,
\label{star-non}
\end{gather}
for any $s,t\in I$ such that $s\leq t$.
\end{lem}

\begin{proof}
The proof is easily done by using the induction on $t-s$.
\end{proof}

The following describes the explicit action of $e_i^c$ on epsilon
systems which is not given in Def\/inition~\ref{def-ep}.
\begin{pro}\label{act-ep}\it
We have the following formula for $s\leq t$:
\begin{gather}
\vep_{[s,t]}(e^c_{t+1}(x))=\vep_{[s,t]}(x)
+\frac{(c-1)\vep_{[s,t+1]}(x)}{\vep_{t+1}(x)},
\label{act-non1}\\
\vep_{[s,t]}(e^c_{s-1}(x))=c\vep_{[s,t]}(x)
+\frac{(1-c)\vep_{[s-1,t]}(x)}{\vep_{s-1}(x)},
\label{act-non2}\\
\vep^*_{[s,t]}(e^c_{t+1}(x))=c\vep^*_{[s,t]}(x)
+\frac{(1-c)\vep^*_{[s,t+1]}(x)}{\vep_{t+1}(x)},
\label{act-*1}\\
\vep^*_{[s,t]}(e^c_{s-1}(x))=\vep^*_{[s,t]}(x)
+\frac{(c-1)\vep^*_{[s-1,t]}(x)}{\vep_{s-1}(x)}.
\label{act-*2}
\end{gather}
\end{pro}

\begin{proof}
Let us only show (\ref{act-non1}) since the others are shown similarly.
It follows from (\ref{non-star}) for $[s,t+1]$ that
\[
 \vep_{[s,t]}(x)=\frac{1}{\vep_{t+1}(x)}\left(\vep_{[s,t+1]}(x)+
\sum_{i=s-1}^{t-1}(-1)^{t-i-1}\vep_{[s,i]}(x)
\vep^*_{[i+1,t+1]}(x)\right).
\]
Then applying $e_t^c$ to $x$, we get
\begin{gather*}
 \vep_{[s,t]}(e^c_{t+1}x) = \frac{1}{c^{-1}\vep_{t+1}(x)}
\left(\vep_{[s,t+1]}(x)+
c^{-1}\sum_{i=s-1}^{t-1}(-1)^{t-i-1}\vep_{[s,i]}(x)
\vep^*_{[i+1,t+1]}(x)\right)\\
\phantom{\vep_{[s,t]}(e^c_{t+1}x)}{} = \frac{1}{\vep_{t+1}(x)}
(c\vep_{[s,t+1]}(x)+\vep_{[s,t]}(x)\vep_{t+1}(x)-\vep_{[s,t+1]}(x))\\
\phantom{\vep_{[s,t]}(e^c_{t+1}x)}{}  = \vep_{[s,t]}(x)
+\frac{(c-1)\vep_{[s,t+1]}(x)}{\vep_{t+1}(x)}.\tag*{\qed}
\end{gather*}
  \renewcommand{\qed}{}
\end{proof}
By this result, we know that all explicit forms of the action by $e_i^c$
on epsilon systems.

\subsection{\bf The proof of Proposition \ref{well-def}}

Let us prove Proposition~\ref{well-def}.
It is trivial to show (\ref{well-0}).
As for (\ref{well-1}), the crucial cases are: for $J=[s,t]$
\begin{gather}
 \vep_J(e^{c_1}_{s-1}e^{c_1c_2}_{s}e^{c_2}_{s-1}x)
=\vep_J(e^{c_2}_{s}e^{c_1c_2}_{s-1}e^{c_1}_{s}x),
\label{well-s}\\
 \vep_J(e^{c_1}_{t}e^{c_1c_2}_{t+1}e^{c_2}_{t}x)
=\vep_J(e^{c_2}_{t+1}e^{c_1c_2}_{t}e^{c_1}_{t+1}x),
\label{well-t}\\
 \vep^*_J(e^{c_1}_{s-1}e^{c_1c_2}_{s}e^{c_2}_{s-1}x)
=\vep^*_J(e^{c_2}_{s}e^{c_1c_2}_{s-1}e^{c_1}_{s}x),
\label{well*-s}\\
 \vep^*_J(e^{c_1}_{t}e^{c_1c_2}_{t+1}e^{c_2}_{t}x)
=\vep^*_J(e^{c_2}_{t+1}e^{c_1c_2}_{t}e^{c_1}_{t+1}x).
\label{well*-t}
\end{gather}
Using the results in Proposition \ref{act-ep},
let us show (\ref{well-s})
\begin{gather*}
 \vep_{[s,t]}(e^{c_1}_{s-1}e^{c_1c_2}_{s}e^{c_2}_{s-1}x)
=c_1\vep_{[s,t]}(e^{c_1c_2}_{s}e^{c_2}_{s-1}x)+
\frac{(1-c_1)\vep_{[s-1,t]}(e^{c_1c_2}_{s}e^{c_2}_{s-1}x)}
{\vep_{s-1}(e^{c_1c_2}_{s}e^{c_2}_{s-1}x)}\\
\quad
 =c_2^{-1}\left(c_2\vep_{[s,t]}(x)+\frac{(1-c_2)\vep_{[s-1,t]}(x)}
{\vep_{s-1}(x)}\right)+
\frac{(1-c_1)\vep_{[s-1,t]}(x)(\vep_{[s-1,s]}(x)+c_2\vep_{[s-1,s]}^*(x))}
{c_2\vep_{s-1}(x)(c_1\vep_{[s-1,s]}(x)+\vep_{[s-1,s]}^*(x))}\\
\quad
 =\vep_{[s,t]}(x)+\frac{\vep_{[s-1,t]}(x)}{c_2\vep_{s-1}(x)}
\left(1-c_2+
\frac{(1-c_1)(\vep_{[s-1,s]}(x)+c_2\vep_{[s-1,s]}^*(x))}
{c_1\vep_{[s-1,s]}(x)+\vep_{[s-1,s]}^*(x)}
\right)\\
\quad
 =\vep_{[s,t]}(x)+\frac{\vep_{[s-1,t]}(x)}{c_2\vep_{s-1}(x)}
\frac{(1-c_1c_2)(\vep_{[s-1,s]}(x)+\vep^*_{[s-1,s]}(x))}
{c_1\vep_{[s-1,s]}(x)+\vep^*_{[s-1,s]}(x))}\\
\quad
 =\vep_{[s,t]}(x)+\frac{(1-c_1c_2)\vep_{[s-1,t]}(x)\vep_s(x)}
{c_2(c_1\vep_{[s-1,s]}(x)+\vep^*_{[s-1,s]}(x))},
\end{gather*}
where the last equality is derived from
$\vep_{[s-1,s]}(x)+\vep^*_{[s-1,s]}(x)=\vep_{s-1}(x)\vep_{s}(x)$.
We also have
\begin{gather*}
 \vep_{[s,t]}(e^{c_2}_{s}e^{c_1c_2}_{s-1}e^{c_1}_{s}x)
=c_2^{-1}\vep_{[s,t]}(e^{c_1c_2}_{s-1}e^{c_1}_{s}x)
=c_2^{-1}\left(c_1c_2\vep_{[s,t]}(e_s^{c_1}x)+
\frac{(1-c_1c_2)\vep_{[s-1,t]}(e_s^{c_1}x)}
{\vep_{s-1}(e_s^{c_1}x)}\right)\\
 \phantom{\vep_{[s,t]}(e^{c_2}_{s}e^{c_1c_2}_{s-1}e^{c_1}_{s}x)}{}
 =\vep_{[s,t]}(x)+\frac{(1-c_1c_2)\vep_{[s-1,t]}(x)\vep_s(x)}
{c_2(c_1\vep_{[s-1,s]}(x)+\vep^*_{[s-1,s]}(x))}.
\end{gather*}
Thus, we obtained (\ref{well-s}). The others are also shown
by direct calculations:
\begin{gather*}
 \vep_{[s,t]}(e^{c_1}_{t}e^{c_1c_2}_{t+1}e^{c_2}_{t}x)
=\vep_{[s,t]}(x)+\frac{(c_1c_2-1)\vep_{[s,t+1]}(x)\vep_t(x)}
{\vep_{[t,t+1]}(x)+c_2\vep^*_{[t,t+1]}(x)}
=\vep_{[s,t]}(e^{c_2}_{t+1}e^{c_1c_2}_{t}e^{c_1}_{t+1}x),
\\
 \vep^*_{[s,t]}(e^{c_1}_{s-1}e^{c_1c_2}_{s}e^{c_2}_{s-1}x)
=\vep^*_{[s,t]}(x)+\frac{(c_1c_2-1)\vep^*_{[s-1,t]}(x)\vep_s(x)}
{c_1\vep_{[s-1,s]}(x)+\vep^*_{[s-1,s]}(x)}
=\vep^*_{[s,t]}(e^{c_2}_{s}e^{c_1c_2}_{s-1}e^{c_1}_{s}x),
\\
 \vep^*_{[s,t]}(e^{c_1}_{t}e^{c_1c_2}_{t+1}e^{c_2}_{t}x)
=\vep^*_{[s,t]}(x)+\frac{(1-c_1c_2)\vep^*_{[s,t+1]}(x)\vep_t(x)}
{c_1(\vep_{[t,t+1]}(x)+c_2\vep^*_{[t,t+1]}(x))}
=\vep^*_{[s,t]}(e^{c_2}_{t+1}e^{c_1c_2}_{t}e^{c_1}_{t+1}x).\tag*{\qed}
\end{gather*}
Let $\ge$ be a Kac--Moody Lie algebra associated with the index set
$I$ and $\ge_J$ be a subalgebra associated with a subset $J\subset I$.
Let ${\bf X}=(X,\{\gamma_i\},\{\vep_i\},\{e_i\})_{i\in I}$ be a $\ge$-geometric
crystal. Then it has naturally a $\ge_J$-geometric crystal structure and
denote it by ${\bf X}_{J}$.
\begin{df}
In the above setting, if $\ge_J$ is isomorphic to the Lie algebra of
type $A_n$ for some~$n$ and
the geometric crystal ${\bf X}_J$ has an epsilon system
$E_{{\bf X}_J}$  of type $A_n$,
then we call it a~{\it local epsilon system} of type $A_n$ associated with the
index set $J$.
\end{df}

\noindent
{\bf Remark.} An $\vep$-geometric crystal has naturally a local epsilon
system associated with each sub-interval of $I$.
\begin{ex}
In Example \ref{ex-1}(ii),
$\{\vep_2,\vep_3,\vep_{23},\vep_{23}^*\}\subset E$ is a local epsilon
system of type $A_2$ associated with the interval $\{2,3\}$.
\end{ex}

\subsection{Product structures on epsilon systems}
\label{subsec-prod}

\begin{thm}\label{thm-prod}\it
Let $\bbX$ and $\bbY$ be $\vep$-geometric crystals.
Suppose that the product
$\bbX\times \bbY$ has a~geometric crystal structure. Then
$\bbX\times \bbY$ turns out to be an $\vep$-geometric crystal
as follows: for $\vep$-systems
$E_\bbX:=\{\vep^X_J,\vep^{X*}_J\}_{J\in\cal J}$ of
$\bbX$ and
$E_\bbY:=\{\vep^Y_J,\vep^{Y*}_J\}_{J\in\cal J}$ of $\bbY$, set
\begin{gather}
\vep_{[s,t]}(x,y):=\sum_{k=s-1}^{t}\frac{\vep^Y_{[s,k]}(y)
\vep^X_{[k+1,t]}(x)}{\prod\limits_{j=s}^{k}\gamma^X_j(x)},
\label{prod-xy}\\
\vep^*_{[s,t]}(x,y):=
\sum_{k=s-1}^{t}
\frac{\vep^{*\,X}_{[s,k]}(x)
\vep^{*\,Y}_{[k+1,t]}(y)}{\prod\limits_{j=k+1}^{t} \gamma^X_{j}(x)}.
\label{prod-xy*}
\end{gather}
Then $E_{\bbX\times \bbY}:=\{\vep_{[s,t]}(x,y),
\vep^*_{[s,t]}(x,y)\}_{[s,t]\in\cJ}$ defines an epsilon
system of $\bbX\times\bbY$.
\end{thm}

\begin{ex} We have
\[
 \vep_{123}(x,y)=\vep_{123}(x)
+\frac{\vep_1(y)\vep_{23}(x)}{\gamma_1(x)}
+\frac{\vep_{12}(y)\vep_3(x)}{\gamma_1(x)\gamma_2(x)}
+\frac{\vep_{123}(y)}{\gamma_1(x)\gamma_2(x)\gamma_3(x)}.
\]
\end{ex}
\begin{proof}
First, we shall show (\ref{ep-ei}).
For $c\in\bbC^\times$ and $(x,y)\in X\times Y$, set
$e_i^c(x,y)=(e_i^{c_1}x,e_i^{c_2}y)$ where
\[
 c_1:=\frac{c\vp_i(x)+\vep_i(y)}{\vp_i(x)+\vep_i(y)},\qquad
c_2:=\frac{c}{c_1}, \qquad \vp_i(x):=\vep_i(x)\gamma_i(x).
\]
Let us see $\vep_{[s,t]}(e_s^c(x,y))=c^{-1}\vep_{[s,t]}(x,y)$.
Each summand $\frac{\vep^Y_{[s,k]}(y)
\vep^X_{[k+1,t]}(x)}{\prod\limits_{j=s}^{k}\gamma^X_j(x)}$
in (\ref{prod-xy}) is changed by the action of $e_s^c$ as
follows (we omit the superscripts $\bbX$ and $\bbY$):
\begin{gather}
\vep_{[s,t]}(e_s^{c_1}x)=c_1^{-1}\vep_{[s,t]}(x),\qquad k=s-1,
\label{k0}\\
\frac{\vep_s(e_s^{c_2}y)\vep_{[s+1,t]}(e_s^{c_1}x)}
{\gamma_s(e_s^{c_1}x)}=
\frac{\vep_s(y)}{c_1^2c_2\gamma_s(x)}
\left(c_1\vep_{[s+1,t]}(x)+\frac{(1-c_1)\vep_{[s,t]}(x)}{\vep_{s(x)}}\right),
\qquad k=s,
\label{k1}\\
 \frac{\vep_{[s,k]}(e_s^{c_2}y)\vep_{[k+1,t]}(e_s^{c_1}x)}
{\gamma_s(e_s^{c_1}x)\cd\gamma_k(e_s^{c_1}x)}
=\frac{\vep_{[s,k]}(y)\vep_{[k+1,t]}(x)}
{c\gamma_s(x)\cd\gamma_k(x)},\qquad k>s,\nonumber
\end{gather}
where the second formula is obtained by (\ref{act-non2}).
Now, taking the summation of (\ref{k0}) and (\ref{k1}), we have
\begin{gather*}
 \vep_{[s,t]}(x)\left(c_1^{-1}+\frac{(1-c_1)\vep_s(y)}
{cc_1\vp_s(x)}\right)+\frac{\vep_s(y)\vep_{[s+1,t]}(x)}
{c\gamma_s(x)}\\
\qquad {}=
c_1^{-1}\vep_{[s,t]}(x)\frac{c\vp_s(x)+\vep_s(y)}
{c(\vp_s(x)+\vep_s(y))}+\frac{\vep_s(y)\vep_{[s+1,t]}(x)}
{c\gamma_s(x)}
=c^{-1}\left(\vep_{[s,t]}(x)+
\frac{\vep_s(y)\vep_{[s+1,t]}(x)}
{\gamma_s(x)}\right).
\end{gather*}
Thus we have $\vep_{[s,t]}(e_s^c(x,y))=c^{-1}\vep_{[s,t]}(x)$.
The others are obtained by the similar argument.

Next, let us show (\ref{ep*ep}).
Set the right hand-side of (\ref{prod-xy}) (resp.\ (\ref{prod-xy*}))
$X_{[s,t]}(x,y)$ (resp.\ $X^*_{[s,t]}(x,y)$).
By Lemma \ref{equiv}, it suf\/f\/ices to show that the relation
(\ref{non-star}) holds
for $X_{[s,t]}(x,y)$ and~$X^*_{[s,t]}(x,y)$
\begin{gather*}
 \sum_{j=s-1}^t(-1)^jX_{[s,j]}(x,y)X^*_{[j+1,t]}(x,y) \\
\qquad {} =\sum_{j=s-1}^t(-1)^j\left(\sum_{k=s-1}^{j}\frac{\vep_{[s,k]}(y)
\vep_{[k+1,j]}(x)}{\prod\limits_{p=s}^{k}\gamma_p(x)}\right)
\left(\sum_{m=j}^{t}
\frac{\vep^{*}_{[j+1,m]}(x)
\vep^{*}_{[m+1,t]}(y)}{\prod\limits_{q=m+1}^{t} \gamma_{q}(x)}\right) \\
\qquad{} =\sum_{s-1\leq k\leq m\leq t}
\frac{\vep_{[s,k]}(y)}{\prod\limits_{p=s}^{k}\gamma_p(x)}
\frac{\vep^{*}_{[m+1,t]}(y)}{\prod\limits_{q=m+1}^{t} \gamma_{q}(x)}
\left(\sum_{j=k}^m(-1)^j\vep_{[k+1,j]}(x)
\vep^{*}_{[j+1,m]}(x)\right)\\
\qquad{}=\frac{1}{{\prod\limits_{p=s}^{t}\gamma_p(x)}}
\sum_{s-1\leq k\leq t}
(-1)^k\vep_{[s,k]}(y)\vep^{*}_{[k+1,t]}(y)=0
\end{gather*}
since $\sum\limits_{j=k}^m(-1)^j\vep_{[k+1,j]}(x)
\vep^{*}_{[j+1,m]}(x)=0$ for $k<m$ and
$\sum\limits_{s-1\leq k\leq t}
(-1)^k\vep_{[s,k]}(y)\vep^{*}_{[k+1,t]}(y)=0$
by~(\ref{non-star}).
Hence, we showed that
$\{X_{[s,t]}(x,y),\,
X^*_{[s,t]}(x,y)\, |\, s,t\in I\ (s\leq t)\}$  satisfy (\ref{non-star}),
which is equivalent to (\ref{ep*ep}) by Lemma~\ref{equiv}.
\end{proof}

\begin{ex}
$D^{(1)}_5$-case \cite{KNO2,KOTY}.
Set the index set $I=\{0,1,2,3,4,5\}$. The geometric crystal~$\cB_L$
is def\/ined as follows:
\begin{gather*}
\cB_L:=\{l=(l_1,\dots,l_5,\ovl l_4,\dots,\ovl l_1)\in(\bbC^\times)^9
|l_1l_2\cdots\ovl l_2 \ovl l_1=L\},\\
\vep_0(l)=l_1\left(\frac{l_2}{\ovl l_2}+1\right),\qquad\!\!
\vep_i(l)=\ovl l_i\left(\frac{l_{i+1}}{\ovl l_{i+1}}+1\right)
\quad (i=1,2,3),\qquad\!\!
\vep_4(l)=l_5\ovl l_4,\qquad\!\! \vep_5(l)=\ovl l_4,\!\\
\gamma_0(l)=\frac{\ovl l_1\ovl l_2}{l_1l_2},\qquad
\gamma_i(l)=\frac{l_i\ovl l_{i+1}}{\ovl l_il_{i+1}}\quad (i=1,2,3),\qquad
\gamma_4(l)=\frac{l_4}{l_5\ovl l_4},\qquad
\gamma_5(l)=\frac{l_4l_5}{\ovl l_4},\\
e_0^c(l)=(\xi_1^{-1}l_1,c^{-1}\xi_1l_2,\dots,
\xi_1\ovl l_2,c\xi^{-1}_1\ovl l_1),\\
e_i^c(l)=(\dots,c\xi_i^{-1}l_i,c^{-1}\xi_il_{i+1},\dots,
\xi_i\ovl l_{i+1},\xi_i^{-1}\ovl l_i,\dots)\quad (i=1,2,3),\\
e_{4}^c(l)=(\dots,c\cdot l_4,c^{-1}l_5,\dots),\qquad
e_5^c(l)=(\dots,c\cdot l_5,c^{-1}\ovl l_4,\dots).
\end{gather*}
where $\xi_i=\frac{l_{i+1}+c\ovl l_{i+1}}{l_{i+1}+\ovl l_{i+1}}$.
There are several local epsilon systems of type $A_4$
in the $D^{(1)}_5$-geometric crystal $\cB_L$ associated with
e.g., $\{0,2,3,4\},\{0,2,3,5\}\subset I$. Then, we have
\begin{gather*}
\vep_{02}(l)=l_1l_2\left(\frac{l_3}{\ovl l_3}+1\right),\qquad
\vep^*_{02}(l)=l_1\ovl l_2\left(\frac{l_3}{\ovl l_3}+1\right),\qquad
\vep_{i\,i+1}(l)=\ovl l_i l_{i+1}
\left(\frac{l_{i+2}}{\ovl l_{i+2}}+1\right),\\
\vep^*_{i\,i+1}(l)=\ovl l_i\ovl l_{i+1}
\left(\frac{l_{i+2}}{\ovl l_{i+2}}+1\right) \quad (i=1,2),\qquad
\vep_{34}(l)=\ovl l_3l_4l_5,\\
\vep^*_{34}(l)=\ovl l_3\ovl l_4l_5,
\vep_{35}(l)=\ovl l_{3}l_4,\qquad
\vep^*_{35}(l)=\ovl l_{3}\ovl l_4,
\\
 \vep_{023}(l)=l_1l_2l_3\left(\frac{l_4}{\ovl l_4}+1\right),\qquad
\vep^*_{023}(l)=l_1\ovl l_2\ovl l_3
\left(\frac{l_4}{\ovl l_4}+1\right),\\
 \vep_{234}(l)=
\ovl l_2l_3l_4l_5,\qquad
\vep^*_{234}(l)=\ovl l_2\ovl l_3\ovl l_4l_5,\\
 \vep_{0234}(l)=l_1l_2l_3l_4l_5,\qquad
\vep^*_{0234}(l)=l_1\ovl l_2\ovl l_3\ovl l_4l_5, \\
 \vep_{0235}(l)=l_1l_2l_3l_4,\qquad
\vep^*_{0235}(l)=l_1\ovl l_2\ovl l_3\ovl l_4,\quad \dots \quad {\rm etc.}\qquad
\raisebox{22mm}[0pt][0pt]{\SelectTips{cm}{}
\qq\q
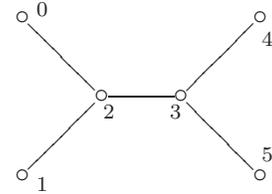
\xymatrix{
*{\circ}<3pt> \ar@{-}[dr]^<{0} &&&*{\circ}<3pt> \ar@{-}[dl]^<{4}\\
&*{\circ}<3pt> \ar@{-}[r]_<{2}
& *{\circ}<3pt>\ar@{}_<{3\,\,\,} &\\
*{\circ}<3pt> \ar@{-}[ur]_<{1}&&&
*{\circ}<3pt> \ar@{-}[ul]_<{5}
}}
\end{gather*}
\end{ex}

\section[Epsilon systems on the Borel subgroup $B^-$]{Epsilon systems on the Borel subgroup $\boldsymbol{B^-}$}
\label{borel}

We shall show that there exists a canonical epsilon system
on the (opposite) Borel subgroup~$B^-$.
Indeed, it is a prototype of general epsilon systems and
derives an epsilon system to
any $A_n$-geometric crystals induced from unipotent crystals.

In this section we shall identify $B^-$ with the set of
lower triangular matrices in $SL_{n+1}(\bbC)$.
Let $T:=\{\text{diag}(t_1,\dots,t_{n+1})\in SL_{n+1}(\bbC)\}$
be the maximal torus in $B^-$.
For $x\in B^-$ there exist a unique unipotent matrix $x_-\in U^-$ and
a unique diagonal matrix $x_0\in T$ such that $x=x_-x_0$.
Then the geometric crystal structure on $B^-$ is described as follows:
For $x\in B^-$, write $x_0={\rm diag}(t_1,\dots,t_{n+1})$ and
$(i,j)$-entry of $x_-$ as
\[
 (x_-)_{i,j}=\begin{cases}u_j&\text{if} \ i=j+1,\\
u_{j,i-1}&\text{if} \ i>j+1,\\
1&\text{if} \ i=j,\\
0&\text{otherwise.}
\end{cases}
\]
For example, for $n=2$-case,
\[
x_-= \begin{pmatrix}1&0&0\\
u_{1}&1&0\\ u_{12}&u_{2}&1
\end{pmatrix},\qq
x_0=\begin{pmatrix}t_1&0&0\\
0&t_2&0\\0&0&t_3
\end{pmatrix},\qq
x= \begin{pmatrix}t_1&0&0\\
t_1u_{1}&t_2&0\\ t_1u_{12}&t_2u_{2}&t_3
\end{pmatrix}.
\]
Now, the rational functions $\gamma_i$ and
$\vep_i$ are given by $\gamma_i(x)=t_i/t_{i+1}$ and $\vep_i(x)=u_i$.
The action $e_i^c$ $(c\in\bbC^\times)$ is given by
\[
 e_i^c(x)=x_i\left(\frac{c-1}{\vep_i(x)}\right)(x)
=x_i\left(\frac{c-1}{\vep_i(x)}\right)\cdot x\cdot
x_i\left(\frac{c^{-1}-1}{\vp_i(x)}\right),
\]
where $x_i(z)={\rm Id}_{n+1}+z E_{i,i+1}\in B\subset SL_{n+1}(\bbC)$.
\begin{pro}\it
In the above setting, the epsilon system on $B^-\subset SL_{n+1}(\bbC)$
are given by:
\begin{equation*}
\vep_{[s,t]}(x)=u_{s,t}, \qquad
\vep^*_{[s,t]}(x)={\rm det}(M_{s,t}),\qquad s<t,
%\label{b-epsilon}
\end{equation*}
where
$M_{s,t}$ is the minor of $x_-$ as:
\[
M_{s,t}:= \begin{pmatrix}
u_s&1&0&\cd&0\\
u_{s,s+1}&u_{s+1}&1&\cd&0\\
&\cd&\cd&\cd&\cd\\
u_{s,t-1}&u_{s+1,t-1}&\cd&\cd&1\\
u_{s,t}&u_{s+1,t}&\cd&u_{t-1,t}&u_t
\end{pmatrix}.
\]
\end{pro}
\begin{proof}
By direct calculations, for $x\in B^-$ we have
\begin{gather*}(e_i^c(x))_- \!=\!\!
\left(\!
\begin{array}{@{}c@{\,}c@{\,}c@{\,}c@{\,}c@{\,}c@{\,}c@{\,}c@{}}
1      &0   &\cd  &\cd  &\cd  &&& \\
u_1    &1   &0    &\cd  &\cd  &&& \\
u_{12} &u_2 &1    &0    &\cd  &&& \\
\cd    &\cd &\cd  &\cd  &\cd  &&\\
u_{1,i-1}{+}\frac{(c-1)u_{1,i}}{u_i}&\cd
%u_{2,i-1}+\frac{(c-1)u_{2,i}}{u_i}
&\cd&
u_{i-1}+\frac{(c-1)u_{i-1,i}}{u_i}&1&\cd&\cd&\cd \\
u_{1,i}  &u_{2,i}  &\cd   & u_{i-1,i} &\frac{u_i}{c} &1 &0&\cd\\
\cd&\cd&\cd&\cd  & \frac{u_{i,i+1}}{c}
&c\left(u_{i+1}{+}\frac{(c^{-1}-1)u_{i,i+1}}{u_i}\right)&1&\cd\\
\cd&\cd&\cd&\cd&\cd&\cd&\cd&\\
u_{1,n}&u_{2,n}&\cd&\cd&\frac{u_{i,n}}{c}
&c\left(u_{i+1,n}{+}\frac{(c^{-1}-1)u_{i,n}}{u_i}\right)&u_n&1
\end{array}\!\right)\!.
\end{gather*}
The formula (\ref{ep-ei}) in Def\/inition \ref{def-ep} is shown by
using this.

We can easily see
\[
 \det(M_{s,t})=u_s\det(M_{s+1,t})-u_{s,s+1}\det(M_{s+2,t})+\cd+
(-1)^{t-s}u_{s,t}.
\]
Thus, by the induction on $t-s$, we obtain (\ref{ep*ep}).
\end{proof}

Next, we shall see an epsilon system on a geometric crystal induced from
a unipotent crystal.
Let $(X,f)$ be a unipotent $SL_{n+1}(\bbC)$ crystal,
where $f:X\to B^-$ is a $U$-morphism.
We assume that any rational function $\vep_i$ is not identically zero.
By Theorem \ref{U-G}, we get the geometric crystal $\bbX=(X,\{e_i^c\},
\{\gamma_i\},\{\vep_i\})$.
\begin{thm}\it
The geometric crystal $\bbX$ as above is an $\vep$-geometric crystal.
Indeed, by setting
\begin{gather}
 \vep^\bbX_{[s,t]}(x):=\vep^{B^-}_{[s,t]}(f(x)),
\label{uni1}\\
 \vep^{*\,\bbX}_{[s,t]}(x)=\sum_{P\in \cP_{[s,t]}}(-1)^{(t-s+1)-l(P)}
\vep^{\bbX}_P(x),
\label{uni2}
\end{gather}
the set $E_\bbX:=\{\vep^\bbX_{[s,t]}(x),\vep^{*\,\bbX}_{[s,t]}(x)|
1\leq s\leq t\leq n\}$ defines an epsilon system on~$\bbX$.
\end{thm}

\begin{proof}
The relation (\ref{ep*ep}) is immediate from (\ref{uni2}).
Since $f(e_i^c(x))=e_i^c(f(x))$ and
$\vep^{B^-}_{[s,t]}$ satisf\/ies~(\ref{ep-ei}),
we can easily obtain (\ref{ep-ei}) for $\vep^\bbX_{[s,t]}$.
\end{proof}

\section{Epsilon systems and tropical R maps}

Let $\bbX=(X,\{e_i^c\},\{\gamma_i\},E_{\bbX})$ and
$\bbY=(Y,\{e_i^c\},\{\gamma_i\},E_{\bbY})$ be $\vep$-geometric crystals,
where $E_\bbX$ (resp.~$E_\bbY$) is the epsilon system of $\bbX$ (resp. $\bbY$).
Suppose that there exists a tropical R-map
\[
 \cR: \ X\times Y\to Y\times X.
\]
Let $E_{\bbX\times \bbY}=\{\vep^{\bbX\times \bbY}_J,
\vep^{*\,\bbX\times \bbY}_J\}_{J\in\cal J}$
(resp.~$E_{\bbY\times \bbX}=\{\vep^{\bbY\times \bbX}_J,
\vep^{*\,\bbY\times \bbX}_J\}_{J\in\cal J}$ ) be the epsilon system on
$\bbX\times\bbY$ (resp.~$\bbY\times\bbX$)
obtained from $E_\bbX$ and $E_\bbY$ by Theorem~\ref{thm-prod}.

\begin{thm}\label{thm-inv}\it
The epsilon systems
$E_{\bbX\times \bbY}$ and
$E_{\bbY\times \bbX}$
are invariant by the action of~$\cR$ in the following sense:
\begin{equation}
\vep^{\bbY\times \bbX}_J(\cR(x,y))=\vep^{\bbX\times \bbY}_J(x,y),\qquad
\vep^{*\,\bbY\times \bbX}_J(\cR(x,y))=\vep^{*\,\bbX\times \bbY}_J(x,y),
\label{formula-inv}
\end{equation}
for any $J\in\cal J$.
\end{thm}
\begin{proof}
Let us show the theorem by the induction on the length of intervals~$J$, denoted $|J|$.
Assume that we have (\ref{formula-inv}) for $|J|<n$.
The case $|J|=1$ is obtained from~(\ref{r2}) in Def\/inition~\ref{def-trop-r}:
the def\/inition of tropical R maps.
Now, we shall show
\begin{equation*}
\vep^{\bbY\times \bbX}_I(\cR(x,y))=
\vep^{\bbX\times \bbY}_I(x,y),\qquad
\vep^{*\,\bbY\times \bbX}_I(\cR(x,y))=
\vep^{*\,\bbX\times \bbY}_I(x,y).
%\label{formula-n}
\end{equation*}
By Theorem \ref{thm-prod}, we have
\begin{equation}
 \vep^{*\,\bbX\times \bbY}_I(x,y)
=\sum_{P\in \cP_I}(-1)^{n-l(P)}
\vep^{\bbX\times \bbY}_P(x,y).
\label{form71}
\end{equation}
In the above summation, if a partition $P\in\cP_I$ is dif\/ferent from
$I=\{1,2,\dots, n\}$, then $\vep^{\bbX\times \bbY}_P(x,y)$ is invariant by
the induction hypothesis.
Thus, we have
\begin{equation}
 \vep^{*\,\bbY\times \bbX}_I(\cR(x,y))
=(-1)^{n-1}\vep^{\bbY\times \bbX}_I(\cR(x,y))+
\sum_{P\in \cP_I,\,\,P\ne I}(-1)^{n-l(P)}
\vep^{\bbX\times \bbY}_P(x,y).
\label{form72}
\end{equation}
It follows from (\ref{form71}) and (\ref{form72}) that
\begin{equation}
\vep^{*\,\bbY\times \bbX}_I(\cR(x,y))-
\vep^{*\,\bbX\times \bbY}_I(x,y)
=(-1)^{n-1}(\vep^{\bbY\times \bbX}_I(\cR(x,y))-
\vep^{\bbX\times \bbY}_I(x,y)).
\label{eq1}
\end{equation}
Here, using (\ref{ep-ei}) we have
\begin{gather*}
 \vep^{*\,\bbY\times \bbX}_I(\cR e_1^c(x,y))
=\vep^{*\,\bbY\times \bbX}_I(e_1^c\cR (x,y))
=\vep^{*\,\bbY\times \bbX}_I(\cR(x,y)),\\
 \vep^{*\,\bbX\times \bbY}_I(e_1^c(x,y))
=\vep^{*\,\bbX\times \bbY}_I(x,y),\\
 \vep^{\bbY\times \bbX}_I(\cR e_1^c(x,y))
=\vep^{\bbY\times \bbX}_I(e_1^c\cR(x,y))
=c^{-1}\vep^{\bbY\times \bbX}_I(\cR(x,y)),\\
 \vep^{\bbX\times \bbY}_I(e_1^c(x,y))=
c^{-1}\vep^{\bbX\times \bbY}_I(x,y).
\end{gather*}
Applying these to (\ref{eq1}), for any $c\in\bbC^{\times}$ one has
\[
 \vep^{*\,\bbY\times \bbX}_I(\cR(x,y))-
\vep^{*\,\bbX\times \bbY}_I(x,y)=(-1)^{n-1}c^{-1}(
\vep^{\bbY\times \bbX}_I(\cR(x,y))-
\vep^{\bbX\times \bbY}_I(x,y)).
\]
Hence, we obtain
\[
  \vep^{*\,\bbY\times \bbX}_I(\cR(x,y))=
\vep^{*\,\bbX\times \bbY}_I(x,y),\qquad
\vep^{\bbY\times \bbX}_I(\cR(x,y))=
\vep^{\bbX\times \bbY}_I(x,y),
\]
which completes the proof.
\end{proof}

Observing this proof, it is easy to get the following:
\begin{cor}\it
Let $\bbX$ $($resp. $\bbY)$ be an $\vep$-geometric crystal with the
epsilon system $E_\bbX=\{\vep_J^\bbX,\vep_J^{*\,\bbX}\}_{J\in\cal J}$
$($resp.\ $E_\bbY=\{\vep_J^\bbY,\vep_J^{*\,\bbY}\}_{J\in\cal J})$
and $F:\bbX\to\bbY$ a homomorphism of geometric
crystals, that is, $F$ is a rational map commuting with the action of
any $e_i$ and preserving the functions~$\vep_i$ and~$\gamma_i$.
Then we obtain for any $J\in\cal J$
\[
\vep^\bbY_J(F(x))= \vep^\bbX_J(x),\qquad
\vep^{*\,\bbY}_J(F(x))= \vep^{*\,\bbX}_J(x).
\]
\end{cor}

\subsection[Application-uniqueness of $A^{(1)}_n$-tropical R map]{Application-uniqueness of $\boldsymbol{A^{(1)}_n}$-tropical R map}

Let  $\{\cB_L\}_{L>0}$ be the  family of $\TY(A,1,n)$-geometric
crystals as in Example \ref{ex-A}.
If we forget the index~0, $\cB_L$ can be seen as an
$A_n$-geometric crystal and is equipped with following local
epsilon system of type~$A_n$:
\begin{gather*}
 \vep_i(l)=l_{i+1}=\vep^*_i(l),\qquad
\vep^*_{[s,t]}(l)=0,\qquad
 \vep_{[s,t]}(l)=l_{s+1}l_{s+2}\cd l_{t+1},
\end{gather*}
for $l=(l_1,\dots,l_{n+1})\in\cB_L$.

\begin{pro}\it
Let $\cR:\cB_L\times\cB_M\to\cB_M\times\cB_L$ be the
tropical R map in Example~{\rm \ref{ex-A}}.
Then~$\cR$ is the unique tropical R maps from
$\cB_L\times\cB_M$ to $\cB_M\times\cB_L$.
\end{pro}
\begin{proof}
For $\cB_L\times \cB_M$, by Theorem \ref{thm-prod} we have
\begin{gather*}
 \vep_i(l,m)=l_{i+1}+\frac{l_{i+1}m_{i+1}}{l_i},\qquad
i=1,\dots,n,\\
 \vep^*_{i,i+1}(l,m)=\vep^*_{i,i+1}(l)+
\frac{\vep_{i+1}(m)\vep_i(l)}{\gamma_{i+1}(l)}+
\frac{\vep^*_{i,i+1}(m)}{\gamma_i(l)\gamma_{i+1}(l)}=l_{i+2}m_{i+2},\qquad
i=1,\dots,n-1.
\end{gather*}

Set $\til l:=L^{\frac{1}{n+1}}$, $\til m:=M^{\frac{1}{n+1}}$,
$l_0:=(\til l,\dots,\til l)\in \cB_L$ and
$m_0:=(\til m,\dots,\til m)\in\cB_M$. Then it is immediate from
the explicit form of $\cR$ as in (\ref{tropr-A}) that
\[
 \cR(l_0,m_0)=(m_0,l_0).
\]
Let $\cR'$ be an arbitrary tropical R-map from
$\cB_L\times\cB_M$ to $\cB_M\times\cB_L$.
Set $(l',m'):=\cR'(l_0,m_0)\in\cB_M\times\cB_L$.
By Theorem \ref{thm-inv}, for any interval $J$ we have
\[
 \vep_J(l',m')=\vep_J(\cR'(l_0,m_0))=\vep_J(l_0,m_0), \qquad
 \vep^*_J(l',m')=\vep^*_J(\cR'(l_0,m_0))=\vep^*_J(l_0,m_0),
\]
and $\gamma_i(l',m')=\gamma_i(l_0,m_0)$.
These equations can be solved uniquely.
Indeed, solving
the system of equations $\vep_i(l',m')=\vep_i(l_0,m_0)$,
$\gamma_i(l',m')=1(=\gamma_i(l_0,m_0))$
$(i=1,\dots,n)$
and $\vep^*_{i,i+1}(l',m')=\vep^*_{i,i+1}(l_0,m_0)$
$(i=1,\dots,n-1)$, that is, for $l'\in\cB_M$ and $m'\in\cB_L$,
\begin{gather*}
  l'_{i+1}+\frac{l'_{i+1}m'_{i+1}}{l'_i}=\til l+\til m,\qquad
l'_im'_i=l'_{i+1}m'_{i+1}, \qquad i=1,\dots,n,\\
 l'_{i+2}m'_{i+2}=\til l\til m, \qquad i=1,\dots,n-1,
\end{gather*}
we obtain the unique solution
$l'_i=\til m$ and $m'_i=\til l$ for any $i=1,\dots,n+1$,
which implies
\[
 (l',m')=(m_0,l_0),
\]
and then $\cR'(l_0,m_0)=\cR(l_0,m_0)$.
According to the remark in Section~\ref{trop-sec}, we have
$\cB_L\cong\cV_L$ and then by Corollary~\ref{preh},
$\cB_L\times\cB_M$ is prehomogeneous.
Therefore,
by Lemma~\ref{uniq} or Theorem~\ref{uniqueness} we have
$\cR=\cR'$.
\end{proof}

\noindent
{\bf Remark.} We expect that this method is applicable to
the tropical R maps of other types~\cite{KNO2}. But we do not have explicit
answers.

\subsection*{Acknowledgments} The author thanks the organizers of the
conference ``Geometric Aspects of Discrete and Ultra-Discrete
Integrable Systems'' for their kind hospitality.
He is supported in part by JSPS Grants-in-Aid for Scientif\/ic
Research $\#$19540050.

\pdfbookmark[1]{References}{ref}
\LastPageEnding

\end{document}